\documentclass{amsart}
\usepackage{amsmath,amssymb, amsfonts}

\newtheorem{propo}{{\bf Proposition}}[section]
\newtheorem{coro}[propo]{{\bf Corollary}}
\newtheorem{remark}[propo]{{\bf Remark}}
\newtheorem{lemma}[propo]{{\bf Lemma}} 
\newtheorem{theor}[propo]{{\bf Theorem}}

\def\N{{\mathbb N}}

\begin{document}
\author{Pilar P\'{a}ez-Guill\'{a}n}
\address{Fakult\"at f\"ur Mathematik, Universit\"at Wien, Oskar-Morgenstern-Platz 1, 1090 Wien, Austria}
\email{maria.pilar.paez.guillan@univie.ac.at}
\author{Salvatore Siciliano}
\address{Dipartimento di Matematica e Fisica ``Ennio De Giorgi", Universit\`{a} del Salento,
Via Provinciale Lecce--Arnesano, 73100--Lecce, Italy}
\email{salvatore.siciliano@unisalento.it}
\author{David A. Towers}
\address{Lancaster University\\
Department of Mathematics and Statistics \\
LA1 4YF Lancaster\\
England}
\email{d.towers@lancaster.ac.uk}


	\begin{abstract}
In this paper we continue the study of the subalgebra lattice of a Leibniz algebra. In particular, we find out that solvable Leibniz algebras with an upper semi-modular lattice are either almost-abelian or have an abelian ideal spanned by the elements with square zero. We also study Leibniz algebras in which every subalgebra is a weak quasi-ideal, as well as modular symmetric Leibniz algebras.
	\end{abstract}	
	
\title[Modularity conditions in Leibniz algebras]{Modularity conditions in Leibniz algebras}
\maketitle

	\section{Introduction}
An algebra $L$ over a field $F$ is called a {\em Leibniz algebra} if, for every $x,y,z \in L$, we have
\[  [x,[y,z]]=[[x,y],z]-[[x,z],y].
\]
In other words, the right multiplication operator $R_x : L \rightarrow L : y\mapsto [y,x]$ is a derivation of $L$. As a result such algebras are sometimes called {\it right} Leibniz algebras, and there is a corresponding notion of {\it left} Leibniz algebras, which satisfy
\[  [x,[y,z]]=[[x,y],z]+[y,[x,z]].
\]
Clearly, the opposite of a right (left) Leibniz algebra is a left (right) Leibniz algebra, so, in most situations, it does not matter which definition we use. Leibniz algebras which satisfy both the right and left identities are sometimes called {\em symmetric} Leibniz algebras.
\par
 
Every Lie algebra is a Leibniz algebra and every Leibniz algebra satisfying $[x,x]=0$ for every element is a Lie algebra. They were introduced in 1965 by Bloh (\cite{bloh}) who called them $D$-algebras, though they attracted more widespread interest, and acquired their current name, through work by Loday and Pirashvili ({\cite{loday1}, \cite{loday2}). They have natural connections to a variety of areas, including algebraic $K$-theory, classical algebraic topology, differential geometry, homological algebra, loop spaces, noncommutative geometry and physics. A number of structural results have been obtained as analogues of corresponding results in Lie algebras.
\par

The {\it Leibniz kernel} is the set $I=$ span$\{x^2:x\in L\}$. Then $I$ is the smallest ideal of $L$ such that $L/I$ is a Lie algebra.
Also $[L,I]=0$.
\par

We define the following series:
\[ L^1=L,\;\; L^{k+1}=[L^k,L] \;\;\; (k\geq 1)\]
and 
\[L^{(0)}=L, \;\; L^{(k+1)}=[L^{(k)},L^{(k)}]\;\;\; (k\geq 0).
\]
Then $L$ is {\em nilpotent of class n} (resp. {\em solvable of derived length n}) if $L^{n+1}=0$ but $L^n\neq 0$ (resp. $ L^{(n)}=0$ but $L^{(n-1)}\neq 0$) for some $n \in \N$. It is straightforward to check that $L$ is nilpotent of class $n$ precisely when every product of $n+1$ elements of $L$ is zero, but some product of $n$ elements is non-zero. Also, $L$ is said to be \emph{supersolvable} if it admits a complete flag made up of ideals of $L$, that is, there exists a chain 
$$
0=L_0\subsetneq L_1 \subsetneq \cdots \subsetneq L_n=L
$$
of ideals of $L$ such that $\dim L_j=j$ for every $0\leq j\leq n$.
\par

The set of subalgebras of a nonassociative algebra forms a lattice under the operations of union, $\cup$, where the union of two subalgebras is the subalgebra generated by their set-theoretic union, and the usual intersection, $\cap$. The relationship between the structure of a Lie algebra $L$ and that
of the lattice ${\mathcal L}(L)$ of all subalgebras of $L$ has been studied by many authors. Much is known about modular subalgebras
(modular elements in ${\mathcal L}(L)$) through a number of investigations including \cite{1,2,3,4,5,6}. Other lattice conditions, together with
their duals, have also been studied. These include semi-modular, upper semi-modular, lower semi-modular, upper modular, lower modular
and their respective duals (see \cite{7} for definitions). For a selection of results on these conditions see \cite{14,8,10,kol,11,sch,13,15,9,12}.
\par

A study of analogous conditions in Leibniz algebras was initiated by Siciliano and Towers in \cite{st}. The subalgebra lattice of a Leibniz algebra is rather different from that of a Lie algebra; for example, in a Lie algebra every element generates a one-dimensional subalgebra, whereas in a Leibniz algebra elements can generate subalgebras of any dimension. So, one would expect different results to hold for Leibniz algebras and this was shown to be the case. The current paper is a further contribution to this study.
\par

Throughout, $L$ will denote a finite-dimensional Leibniz algebra over a field $F$. Let $J=\langle x\in L \mid x^2=0\rangle$. In \cite{st} it was shown that $J$ plays a role in determining when $L$ is modular. We recall that $L$ is said to be \emph{modular} if $\langle U,V\rangle \cap W=\langle U,V \cap W\rangle$ for all subalgebras $U,V,W$ of $L$ with $U\subseteq W$. Moreover, a subalgebra $U$ of $L$ is called {\em upper semi-modular} if it is a maximal subalgebra of $\langle U,B\rangle$ for every subalgebra $B$ of $L$ such that $U\cap B$ is maximal in $B$. We say that $L$ is {\em upper semi-modular} if every subalgebra of $L$ is upper semi-modular in $L$. In Section 2 we make a more detailed study of $J$ for solvable upper semi-modular Leibniz algebras.  We define two types of {\em almost abelian} Leibniz algebras $L=A+Fy$ where $A$ is an abelian ideal and, in each, $[a,y]=a$ for all $a\in A$. We say that $L$ is an {\em almost abelian Lie algebra} if also $[y,a]=-a$, and that it is an {\em almost abelian non-Lie Leibniz algebra} if $[y,a]=0$. It is shown that if $L$ is an upper semi-modular solvable Leibniz algebra then $J$ is an abelian or almost abelian subalgebra of $L$. It is further shown that, in such an algebra $L$, if $J$ is almost abelian, then $L=J$. Finally it is shown that, for such an algebra, if $F$ has characteristic different from two, $J$ is an ideal of $L$. 
\par

The {\em Frattini ideal}, $\phi(L)$, of $L$ is the largest ideal contained in the intersection of all maximal subalgebras of $L$. A subalgebra $U$ of $L$ to be a {\em weak quasi-ideal} of $L$ if $[U,V]+[V,U]\subseteq U+V$ for every subalgebra $V$ of $V$. In Section 3 we study Leibniz algebras $L$ in which every subalgebra is a weak quasi-ideal of $L$. In such an algebra, $L/\phi(L)$ is either abelian or almost abelian. We classify the algebras $L$  in which every subalgebra is a weak quasi-ideal of $L$ for the cases where $L/\phi(L)$ is almost abelian. If $L$ is a non-Lie almost abelian Leibniz algebra, no restrictions on $F$ are needed, but if it is an almost abelian Lie algebra we require $F$ to have characteristic different from two and for $\sqrt{F}\subseteq F$. We close this section by noting that, if $F$ is algebraically closed, then the conditions that $L$ is modular, $L$ is upper semi-modular and $L$ has every subalgebra a weak quasi-ideal are equivalent. 
\par

In the final section we apply the preceding results to classify all symmetric Leibniz algebras over an algebraically closed field of characteristic different from two which are modular.	
\par

We denote by $\langle x,y\rangle$ the subalgebra of $L$ generated by $x,y\in L$. The notation $\oplus$ will denote an algebra direct sum, whilst $\dot{+}$ will denote a direct sum of the underlying vector space structure. Also, $\subseteq$ will denote inclusion and $\subset$ will denote strict inclusion.
	
	\section{Results on $J$}
In this section we study $J$ when $L$ is a solvable upper semi-modular Leibniz algebra. First we need some preliminary results. The first is straightforward to show.

\begin{lemma}\label{kernel} Let $L$ be a Leibniz algebra over any field $F$ with Leibniz kernel $I$ and $\phi(L)\subseteq I$. Then the Leibniz kernel of $L/\phi(L)$ is $I(L/\phi(L))=I/\phi(L)$.
	\end{lemma}

	\begin{lemma}\label{1dim} Let $L$ be a solvable Leibniz algebra in which all proper subalgebras are one-dimensional. Then $L$ is two-dimensional and is a Lie algebra or is cyclic. 
	\end{lemma}
	\begin{proof} If $\langle x\rangle$ is a proper subalgebra of $L$, and so $x^2=0$, for all $x\in L$, then $L$ is a Lie algebra. If $L$ is not abelian, we have that $\dim L^2=1$. But now $Fx+L^2$ is a subalgebra for every $x\notin L^2$, so $L=Fx+L^2$ and $L$ is two-dimensional.
		\par
		
		If there exists $x\in L$ such that $\langle x\rangle$ is not a proper subalgebra of $L$, then $L$ is cyclic. Moreover, it follows from \cite[Proposition 6.1]{st} and its proof that $L$ is two-dimensional.
	\end{proof}
	\medskip
	
	\begin{lemma}\label{two} Let $L$ be an upper semi-modular solvable Leibniz algebra over any field $F$. If $L$ is generated by two distinct one-dimensional subalgebras, then $L$ is two-dimensional.
	\end{lemma}
	\begin{proof} Put $L= \langle x,y\rangle$, where $x,y$ are distinct elements of $L$ with $x^2=y^2=0$, and let $S$ be a non-zero proper subalgebra of $L$. Suppose first that $x\in S$, $y\notin S$. Then $Fx\cap Fy=0$ is covered by $Fy$, so $Fx$ is covered by $L$. It follows that $S=Fx$. 
		\par
		
		Suppose now that $x,y\notin S$, and let $s\in S$. If $s^2=0$, then $Fs\cap Fx=0$ is covered by $Fx$ and $Fs$, so $Fx$ and $Fs$ are covered by $\langle x,s\rangle\subseteq L$. Hence $\langle x,s\rangle=L$ and $S=Fs$. If $s^2\neq 0$, then replace $s$ by $s^2$ in the previous sentence to deduce that $\langle x,s^2\rangle=L$ and $S=Fs^2$, which yields a contradiction, since $s=\lambda s^2$ ($\lambda \in F$) whence $s^2=0$. 
		\par
		
		We have shown that all proper non-zero subalgebras of $L$ are one-dimensional. It follows from Lemma \ref{1dim} that $L$ is two-dimensional.
	\end{proof}
	
	\begin{coro} Let $L$ be an upper semi-modular solvable Leibniz algebra over any field $F$. Then $J = sp\{ x\in L \mid x^2=0\}$.
	\end{coro}
	
	\begin{lemma}\label{three}  Let $L$ be a non-abelian upper semi-modular solvable Leibniz algebra over any field $F$. If $L$ is generated by three distinct one-dimensional subalgebras, then $Z(L)=0$.
	\end{lemma}
	\begin{proof} Let $L= \langle x,y,z\rangle$ where $x,y,z$ are distinct elements of $L$ with $x^2=y^2=z^2=0$ and suppose that $Z(L)\neq 0$. By Lemma \ref{two} and without loss of generality we may assume that $x\in Z(L)$ and that $\langle y,z\rangle$ is almost abelian with $[z,y]=z$. Then $Fy \cap F(x+z)=0$ is maximal in $Fy$ and $F(x+z)$, but neither of these is maximal in $\langle y,x+z\rangle=\langle x,y,z \rangle$, a contradiction.
	\end{proof}
\medskip

Now we have the promised characterisation of $J$.
	
	\begin{theor}\label{abalab} Let $L$ be an upper semi-modular solvable Leibniz algebra over any field. Then $J$ is abelian or almost abelian.
	\end{theor}
	\begin{proof} By Lemma \ref{two}, all subalgebras of $J$ generated by two distinct one-dimensional subalgebras are abelian or almost abelian. Suppose that $\langle y,x_1\rangle$ is almost abelian with $[x_1,y]=x_1$ and $[y,x_1]=-x_1$ or $0$. Write $J=\langle y,x_1,\ldots,x_n\rangle$, where $y,x_1,\ldots,x_n$ are linearly independent and $y^2=x_1^2=\cdots=x_n^2=0$. We distinguish two cases depending on the value of $[y,x_1]$.

		\medskip

		\noindent {\bf Case 1:}  Suppose first that $\langle y,x_1\rangle$ is an almost abelian Lie algebra, and suppose too that $[x_i,y]=[y,x_i]=0$ for some $2\leq i\leq n$. By Lemma \ref{three} we must have that $[x_1,x_i]\neq 0$ or $[x_i,x_1]\neq 0$, so $\langle x_1,x
		_i\rangle$ is almost abelian. Now $[x_1,x_i]=\lambda x_1+\mu x_i$ for some $\lambda, \mu\in F$. But $[x_1,x_i]=[[x_1,y],x_i]=[[x_1,x_i],y]=\lambda[x_1,y]+\mu[x_i,y]=\lambda x_1$. Similarly, $[x_i,x_1]=\alpha x_1$. As $[x_1,x_i]+[x_i,x_1]\in I$, we have that $0=[y,[x_1,x_i]]+[y,[x_i,x_1]]=-(\lambda+\alpha)x_1$, whence $\alpha=-\lambda\neq 0$. But now $y+\lambda^{-1}x_i\in Z(\langle y,x_1,x_i\rangle)=0$ by Lemma \ref{three}, a contradiction.  Hence, $[x_i,y]\neq 0$ or $[y,x_i]\neq 0$. Put $[x_i,y]=\gamma x_i+\delta y$, $[y,x_i]=\epsilon x_i+\zeta y$. Then
		\begin{align*} \delta [x_i,y]=&[x_i,[x_i,y]]=-[[x_i,y],x_i]=-\delta [y,x_i], \\
			\gamma [y,x_i]=&[y,[x_i,y]]=[[y,x_i],y]=\epsilon [x_i,y], \\
			\epsilon [y,x_i]=&[y,[y,x_i]]=-[[y,x_i],y]=-\epsilon [x_i,y].
		\end{align*} If $[x_i,y]\neq -[y,x_i]$, then $\delta=\epsilon=0$ and either $\gamma=0$ or $\zeta=0$  (but not both). If $\gamma=0$, then $y\in I$, which contradicts $[x_1,y]=x_1$, so $\zeta=0$. Then $x_i\in I$, so $[x_1,x_i]=0$. Let $[x_i,x_1]=\alpha x_i+\beta x_1$. Then $$\gamma [x_i,x_1]=[[x_i,y],x_1]=[x_i,[y,x_1]]+[[x_i,x_1],y]=-[x_i,x_1]+\alpha[x_i,y]+\beta[x_1,y],$$
	so	$$ \gamma\alpha x_i +\gamma\beta x_1=-\alpha x_i-\beta x_1+\alpha\gamma x_i+\beta x_1,$$ whence $\gamma\beta x_1 =-\alpha x_i$, and $\alpha=\beta=0$. But now $Fy\cap F(x_1+x_i )=0$ is maximal in $Fy$ and $F(x_1+x_i)$, but neither is maximal in $\langle y, x_1+x_i \rangle = \langle y,x_1,x_i \rangle$, a contradiction.
		\par
		
		Hence $[x_i,y]=-[y,x_i]=\gamma x_i + \delta y$. Let $[x_1,x_i]=\eta x_1 +\phi x_i$. Now, $F(y+x_i)\cap Fx_1=0$ is maximal in $F(y+x_i)$ and $Fx_1$. But, $[x_i+y,x_1]=\alpha x_i+\beta x_1 -x_1$ and $[x_1,x_i+y]=\eta  x_1 +\phi x_i+x_1$, so, if $\alpha, \phi \neq 0$, neither $F(y+x_i)$ nor $Fx_1$ is maximal in $\langle y+x_i,x_1\rangle=\langle y, x_1,x_i\rangle$. Hence $[x_i,x_1]=\beta x_1$ and $[x_1,x_i]=\eta x_1$.
		\par
		
		Then 
		$$[[x_i,y],x_1]=[[x_i,[y,x_1]]+[[x_i,x_1],y]=-[x_i,x_1]+\beta[x_1,y]=0, $$ 
		so $0=\gamma[x_i,x_1]+\delta[y,x_1]=(\beta\gamma-\delta)x_1$, whence $\beta\gamma=\delta$. Also 
		$$[x_1,[x_i,y]]=[[x_1,x_i],y]-[[x_1,y],x_i]=\eta[x_1,y]-[x_1,x_i]=0,$$
		so $0=\gamma[x_1,x_i]+\delta[x_1,y]=(\gamma\eta+\delta)x_1$, whence $-\gamma\eta=\delta$. Now $\gamma=0$ implies that $\delta =0$ and $[x,y_i]=[y_i,x]=0$, a contradiction. It follows that $\eta=-\beta$.
		\par
		
		But now $Fy\cap F(x_1+x_i)=0$ is maximal in $Fy$ and $F(x_1+x_i)$ and neither is maximal in $\langle y,x_1+x_i \rangle=\langle y,x_1,x_i \rangle$ unless $\gamma=1$. Thus, $[[x_i,y],y]=[x_i,y]$ and $[[x_i,y],x_1]=[x_i,[y,x_1]]+[[x_i,x_1],y]=-[x_i,x_1]+\beta[x_1,y]=0$, so, by replacing $x_i$ by $[x_i,y]$ we may assume that $[x_i,y]=x_i$ and $[x_i,x_1]=0$. Finally, $[x_i,x_j]=\lambda x_i+\mu x_j= [[x_i,x_j],y]=[x_i,[x_j,y]]+[[x_i,y],x_j] =[x_i,x_j]+[x_i,x_j]$, so $[x_i,x_j]=0$ and $J$ is an almost abelian Lie algebra.
		\medskip
		
		\noindent {\bf Case 2:} Suppose now that $\langle y,x_1\rangle$ is a non-Lie almost abelian Leibniz algebra, so $[x_1,y]=x_1$ and $[y,x_1]=0$. Then $x_1=[x_1,y]+[y,x_1]\in I$, so $[x_i,x_1]=0$ for all $1\leq i\leq n$. Also, $[x_1,x_i]=\eta x_1+\phi x_i\in I$, so $\phi[x_1,x_i]=0$, which yields that $[x_1,x_i]=\eta x_1$. 
		\par
		
		Now, if $[x_i,y]=\gamma x_i+\delta y$, we have $\delta [x_i,y]=[x_i,[x_i,y]]=-[[x_i,y],x_i]=-\delta[y,x_i]$. If $[x_i,y]=-[y,x_i]$ we have that $Fy\cap F(x_1+x_i-\eta y)=0$ is maximal in $Fy$ and $F(x_1+x_i-\eta y)$ and neither is maximal in $\langle y,x_1+x_i-\eta y\rangle$, a contradiction, so $\delta=0$.

		\par
		
		If $\gamma =0$, then $[y,x_i]=\epsilon x_i+\zeta y\in I$, so $\epsilon [y,x_i]=0$. If $[y,x_i]\neq 0$, then $y\in I$, contradicting $[x_1,y]=x_1$, so $[y,x_i]=0$.  
		We thus have $[x_i,y]=0=[y,x_i]$, and we get a contradiction as in the previous paragraph. 
		\par
		
		So, $\gamma\neq 0$, giving $x=(\gamma+\epsilon)x_i+\zeta y\in I$, whence $0=[x_i,x]=\zeta\gamma x_i$. So $\zeta=0$ and either $x_i\in I$ or $\gamma=-\epsilon$. 
		If $\gamma=-\epsilon$, then $[x_i,y]=-[y,x_i]$, a contradiction. Then $x_i\in I$ and consequently $\eta=\epsilon=0$; also, $[x_i,x_j]=0$ for all $1\leq i,j\leq n$. Finally, we can set $\gamma=1$ and $L$ is a non-Lie almost abelian Leibniz algebra.	
	\end{proof}

	\begin{propo}\label{usm2} Let $L$ be an upper semi-modular solvable Leibniz algebra. Then $L/I$ is abelian or almost abelian.
	\end{propo}
	\begin{proof} This follows from \cite[Lemma 4.2]{st} and \cite[Theorem 2.3]{kol}.
	\end{proof}

	\begin{theor}\label{alab}  Let $L$ be an upper semi-modular solvable Leibniz algebra over any field $F$. If $J$ is an almost abelian Leibniz algebra, then $L=J$.
	\end{theor}
	\begin{proof} Let $J=Fy+A$ where $[a,y]=a$  and $[y,a]=-a$ or $0$ for all $a\in A$, and $A$ is abelian. Let $x\in L\setminus J$. Then $x^2\in I\subseteq J$, so $x^2=\lambda y+a$ for some $a\in A$, $\lambda \in F$. Now $0=[a,x^2]=\lambda a$. If $a=0$, then $\lambda y\in I$ and $\lambda[ a',y]=0$ for some $a'\in A$, $a'\neq 0$, which implies that $\lambda=0$ and $x\in J$, a contradiction. Thus, $x^2\in A$ and $I\subseteq A$. Now we distinguish two cases. 
		\par
		
		If $J$ is a Lie algebra, then $0=[y,x^2]=-x^2$, so $x\in J$, a contradiction. Hence, in this case, $L=J$. 
		\par
		
		Thus, suppose that $[y,a]=0$ for all $a\in A$. Then $a=[a,y]\in I$ for all $a\in A$, so $A\subseteq I$, whence $A=I$. Suppose first that $L/I$ is an almost abelian Lie algebra. If $y+I\in L^2/I$, then there exists $ z +I\in L/I$ such that $[y,x ]=y+a'$ for some $a'\in A$. But now, if $a$ is any element of $A$, we have that $a=[a,y]=[a,[y, z]]=[[a,y], z]-[[a, z],y]=[a, z]-[a, z]=0$, contradicting the fact that $J$ is almost abelian. It follows that $[ z,y]= z+a_{ z}$ for some $a_{ z}\in A$, for all $ z+I\in L^2/I$. Note that we can assume that $x\in L^2/I$. Now $[a,x]=[a,[x,y]]=[[a,x],y]-[[a,y],x]=[a,x]-[a,x]=0$ for all $a\in A$; similarly, $[x,a]=0$ for all $a\in A$. Then $x^2=[x^2,y]=[x,[x,y]]+[[x,y],x]=x^2+x^2$ and $x^2=0$. It follows that $x\in J$, a contradiction.
		\par
		
		Finally, suppose that $L/I$ is abelian. In this case, both $[x,y]$ and $[y,x]$ belong to $I$. Then $x^2=[x^2,y]=[x,[x,y]]+[[x,y],x]=[[x,y],x]$;  $[x,[x,y]]=0$  and $[y,x]=[[y,x],y]=[y,[x,y]]+[y^2,x]=0$.  Thus, $F[x,y]\cap F(x-[x,y])$ is maximal in $F(x-[x,y])$, so $F[x,y]$ is maximal in $\langle [x,y],x-[x,y]\rangle=\langle x,[x,y]\rangle$. But $F[x,y]\subseteq \langle x^2,[x,y]\rangle\subset \langle x,[x,y]\rangle$. It follows that $x^2=\alpha [x,y]$ for some $\alpha \in F$. 
		Then $Fy\cap F(x-\alpha y)=0$ is maximal in $F(x-\alpha y)$.
		Hence $Fy$ is maximal in $\langle y,x-\alpha y\rangle=\langle x,y\rangle$. However, $Fy\subset \langle x^2,y\rangle\subset \langle x,y\rangle$ unless $x^2=0$ and $x\in J$, completing the proof.
	\end{proof}
	
	\begin{theor}\label{ideal}  Let $L$ be a solvable upper semi-modular Leibniz algebra over a field $F$ of characteristic different from two. Then $J$ is an ideal of $L$.
	\end{theor}
	\begin{proof} If $J$ is almost abelian, then the result follows from Theorem \ref{alab}, so suppose that $J$ is abelian. If $L/I$ is abelian,then $$[L,J]+[J,L]\subseteq L^2\subseteq I\subseteq J,$$ and $J$ is an ideal of $L$.
		\par
		
		So assume that $L/I$ is almost abelian, say $L/I=A/I+(Fy+I)/I$ where $A^2\subseteq I\subseteq J$ and, for every $a\in A$, $[a,y]=a+i_a$, $[y,a]=-a + i_a'$ for some $i_a,i_a'\in I$. Suppose that $y\in J$. Then, for every $a\in A$, $$0=[a^2,y]=[a,[a,y]]+[[a,y],a]=[a,a+i_a]+[a+i_a,a]=2a^2+[i_a,a].$$ But $[i_a,a]=[i_a,[a,y]]=[[i_a,a],y]-[[i_a,y],a]=0$, since $y\in J$ and $J$ is abelian, so $2a^2=0$ and $a\in J$. It follows that $J=L$ and the result holds in this case.
		\par
		
		Suppose finally that $y\notin J$. Let $j=a+\lambda y+i\in J$, where $0\neq a\in A$, $i\in I$, $\lambda\in F$. Then $[j,a]=a^2+\lambda[y,a]+[i,a]=-\lambda a+i'$, $[a,j]=a^2+\lambda [a,y]=\lambda a+i''$, where $i',i''\in I$. Hence $\lambda[j,a]=[j,[a,j]]=[[j,a],j]=- \lambda[a,j]$, so $i''=- i'$ or $\lambda=0$. Suppose $\lambda\neq 0$. It follows that $\lambda a^2=[[a,j],a]+[i',a]=[a,[j,a]]+[i',a]=-\lambda a^2+[a,i']+[i',a]$. But $[i',a]=\lambda^{-1}[i',[a,j]]=\lambda^{-1}[[i',a],j]=0$ and $[a,i']=-\lambda^{-1}[[j,a],i']=-\lambda^{-1}[j,[a,i']]=0$, so $2\lambda a^2=0$ and $a\in J$. Thus, $\lambda y\in J$, a contradiction, whence $\lambda=0$. Then $[A,j]+[j,A]\subseteq I\subseteq J$, $[y,j]=-a+i_a'\in J$ and $[j,y]=a+i_a+[i,y]\in J$, so $J$ is an ideal of $L$.
		
	\end{proof}

	\section{Algebras with every subalgebra a weak quasi-ideal}
In this section we study Leibniz algebras $L$ in which every subalgebra is a weak quasi-ideal of $L$.

\begin{lemma}\label{qi} Let $L$ be a Leibniz algebra over any field $F$. The following conditions are equivalent:
	\begin{itemize}
		\item[(i)] every subalgebra of $L$ is a weak quasi-ideal of $L$; and
		\item[(ii)] $[x,y]\in \langle x\rangle  +\langle y\rangle  $ for all $x,y\in L$.
	\end{itemize}
\end{lemma}

\begin{lemma} Let $L$ be a Leibniz algebra in which every subalgebra is a weak quasi-ideal. Then $L/\phi(L)$ is abelian or almost abelian.
\end{lemma}
\begin{proof} Clearly, the hypothesis is factor algebra closed, so every subalgebra of $L/\phi(L)$ is a weak quasi-ideal. If $L/\phi(L)$ is a Lie algebra, then every subalgebra of it is a quasi-ideal and so it is abelian or an almost abelian Lie algebra, by \cite[Theorem 3.8]{amayo}. If it is not a Lie algebra, then every subalgebra of it is modular, and hence upper semi-modular, by \cite[Proposition 6.3]{st}. It follows that it is an almost abelian non-Lie Leibniz algebra, by \cite[Proposition 4.4]{st}.
\end{proof}

We devote the next pages to classify the Leibniz algebras $L$ in which every subalgebra is a weak quasi-ideal and such that $L/\phi(L)$ is almost abelian. If $L/\phi(L)$ is abelian, then $L$ is nilpotent, and the nilpotent Leibniz algebras in which every subalgebra is a weak quasi-ideal of $L$ are more difficult to characterise, as was observed in \cite{st}.

 We will need the following result, which was proved as part of \cite[Proposition 6.3]{st}. Note that the assumption of an algebraically closed field is not needed for this part of that result, as is remarked after its proof.
 
             \begin{lemma}\label{cyclic} Let $L$ be a cyclic Leibniz algebra $L$ of dimension $n$ over a field $F$. Then every subalgebra of $L$ is a weak quasi-ideal of $L$ if and only if it is one of the following two types:
\begin{itemize}
\item[(i)] nilpotent, so $L=\langle a\rangle$ where $[a^i,a]=a^{i+1}$ for $1\leq i\leq n-1$, and all other products are zero; or
\item[(ii)] solvable with $L=\langle a\rangle$ where $[a^i,a]=a^{i+1}$ for $1\leq i\leq n-1$, $[a^n,a]=a^n$, and all other products are zero.
\end{itemize}
            \end{lemma}
            \begin{proof} Suppose first that every subalgebra of $L$ a weak quasi-ideal of $L$. Then $L$ is modular by Lemma \ref{qi}. (Note that the algebraically closed field is not needed for this implication.) Hence $L$ is of the form given in (i) or (ii), by \cite[Proposition 6.1]{st}.
\par

Conversely, suppose that $L$ is of the form (i) or (ii). In case (i), $\phi(L)=L^2$, and in case (ii), $\phi(L)=\sum_{i=2}^n F(a^i-a^{i-1})$ (see the proof of  \cite[Proposition 6.1]{st}). It follows that if $x=\sum_{i=1}^n \lambda_ia^i\in L$ and $\lambda_1\neq 0$, then $L=\langle x\rangle +\phi(L)$, so $\langle x\rangle =L$; if $\lambda_1=0$, then $\langle x \rangle=Fx$.
\par

Let $y=\sum_{i=1}^n \mu_ia^i\in L$. Then $[x,y]=\sum_{i=2}^n\mu_1\lambda_i a^{i+1}\in L$. If $\mu_1\neq 0$ then $[x,y]\in L=\langle y\rangle$; if $\mu_1=0$ then $[x,y]=0$. It follows from  Lemma \ref{qi} that every subalgebra of $L$ is a weak quasi-ideal of $L$.
           \end{proof}

	\begin{lemma}\label{int} Let $L$ be a Leibniz algebra over any field in which every subalgebra is a weak quasi-ideal. If $\phi(L)\neq 0$ then $I\cap \phi(L)\neq 0$.
	\end{lemma}
	\begin{proof} We prove the contrapositive.  Suppose that $I\cap \phi(L)=0$. Then there is a subalgebra $B$ of $L$ such that $L=I\dot{+} B$, by \cite[Lemma 7.2]{frat}. Also, $B$ is an abelian or almost abelian Lie algebra, by \cite[Theorem 3.8]{amayo}. It follows that $\phi(L) \subseteq I$, whence $I\cap \phi(L)=\phi(L)= 0$. 
	\end{proof}
	
	\begin{theor}\label{non-Lie} Let $L$ be a Leibniz algebra over any field $F$. Then every subalgebra is a weak quasi-ideal and $L/\phi(L)$ is a non-Lie almost abelian Leibniz algebra if and only if either
		\begin{itemize}
			\item[(i)] $L$ is a non-Lie almost abelian Leibniz algebra, or
			\item[(ii)] $L=C\dot{+}A$, where $C$ is a cyclic Leibniz algebra with basis $x, \ldots x^k$, $x^{k+1}=x^k$ with $k\geq 2$, $A$ is abelian (possibly $0$) with $[a,x]=a$ for all $a\in A$, all other products zero. 
		\end{itemize}
	\end{theor}
	\begin{proof} Let us prove the necessity. 		
		Suppose first that  $L/\phi(L)$ is cyclic.  Then so is $L$. Moreover, $L$ must be of the form given in \cite[Proposition 6.1]{st}, using Lemma \ref{qi}. But it cannot be nilpotent as that would imply that  $L/\phi(L)$  is nilpotent. It follows that $L$ is of the form (ii) with $A=0$. 
		
		So suppose that $L/\phi(L)$ is not cyclic. 
		Note that $L$ must be supersolvable, since $L/\phi(L)$ is supersolvable. We use induction on $\dim \phi(L)$. The result clearly holds if $\phi(L)=0$, so suppose it holds whenever $\dim \phi(L)\leq m$ ($m\geq 0$) and let $L$ be such that $\dim \phi(L)=m+1$. Let $Fa_n$ be a minimal ideal of $L$ with $Fa_n\subseteq I\cap \phi(L)$, which exists by Lemma \ref{int}.  Then $\phi(L/Fa_n)=\phi(L)/Fa_n$, by \cite[Proposition 4.3]{frat}, so $\dim \phi(L/Fa_n)= m $ and $L/Fa_n$ is of the given form.
		
		\par
		
		Suppose first that it is of the form (i). Then $L$ will have a basis $a_1,\ldots,  a_n $, $x$ where $[a_i,a_j]=[x,a_i]=0$  for all $1\leq i,j\leq n$, by Lemma \ref{kernel}. Also, $[a_i,x]=a_i+\alpha_i a_n$ for $1\leq i\leq n-1$, $[a_n,x]=\alpha_n a_n$  and $x^2=\beta a_n$. 
		\par
		
		Now there is a subalgebra $U$ which is minimal with respect to $L=L^2+U$. It follows from \cite[Lemma 7.1]{frat} that $L^2\cap U\subseteq \phi(U)$,  so $U^2\subseteq L^2\cap U\subseteq \phi(U)$, whence $U$ is nilpotent, by \cite[Theorem 5.3]{barnes}.  As $L/L^2\cong U/L^2\cap U$, we can take $x$ to be nilpotent and two possibilities arise:
		\begin{itemize}
			\item[(a)] $x^2=0$, or
			\item[(b)] we can take $a_n=x^2$ and $[a_n,x]=0$. 
		\end{itemize}
		\par
		
		Suppose first that (a) holds. Then $[a_i,x]= \alpha a_i+\beta x$, which implies that $\alpha=1$ and $\alpha_i=0$ for $1\leq i\leq n-1$. Also, we have that $[a_1+a_n,x]=a_1+\alpha_n a_n=\alpha(a_1+a_n)+\beta x$, whence $\alpha=\alpha_n=1$, $\beta=0$ and $L$ is as in (i).
		\par
		
		So, suppose next that (b) holds. Then $[[a_i,x],x]=[a_i,x]$ for $1\leq i\leq n-1$, so, by replacing $a_i$ by $[a_i,x]$ we can assume that $\alpha_i=0$.  Suppose that $n\geq 3$. Then $$[a_1+x^2, a_2+x]=a_1=\alpha(a_1+x^2)+\beta(a_2+x)+\gamma(a_2+x^2),$$ and there are no values of $\alpha, \beta, \gamma$ satisfying this. If $n=2$ then $L$ is cyclic, generated by $a_1+x$ so this case was already dealt with. If $n=1$ then $L$ is two-dimensional and is nilpotent, so $L/\phi(L)$ is abelian, contradicting the hypothesis.
		\par
		
		Suppose now that $L/Fa_n$ is of the form (ii).
		 Then $L$ will have a basis $a_1,\ldots, a_n$, $x, \ldots x^k$ where $[a_i,x]=a_i+\lambda_i a_n$ for $1\leq i\leq n-1$, $x^{k+1}=x^k+\gamma a_n$ and $[a_n,x]=\delta a_n$. Moreover, $[a_i,a_j] =[x^l,a_i]=0$ for all $1\leq i,j\leq n$ and $1\leq l\leq k$, by Lemma \ref{kernel}. 
		\par

		If $\gamma=0$ then $[a_i,x]=a_i+\lambda_ia_n\in Fa_i+\langle x\rangle$ implies that $\lambda_i=0$ for $1\leq i\leq n-1$. Similarly, $[a_1+a_n,x]=a_1+\delta a_n$ implies that $\delta=1$. Hence, putting $C= Fx+\cdots+Fx^n$ and $A=\sum_{i=1}^nFa_i$ shows that $L$ is as described in (ii).
		\par
		
		If $\gamma \neq 0$ then $a_n\in \langle x \rangle$ and $C=Fx+\cdots +Fx^k+Fa_n$ is a cyclic subalgebra of $L$.  As $C/Fa_n$ is not nilpotent, neither is $C$,  so $C$ must be of the form given in (ii)  by \cite[Proposition 6.1]{st}. Replacing $a_n$ by $x^k+\gamma a_n$ we can assume that $a_n=x^{k+1}$.  Note that $\delta\neq 0$ and we can assume $\delta=1$.  But now  $[a_i,x]=a_i+\lambda_i x^{k+1}$ for $1\leq i\leq n-1$. We claim that $\lambda_i=0$ for all $1\leq i \leq k$. Suppose first that $n=1$ and $\lambda_1\neq 0$. Consider $y=a_1+x$. Then 
		$y^r=a_1+x^r+(r-1)\lambda_ix^k$ for all $r\geq 1$ and consequently $y$ generates the whole algebra. Then $L$ is cyclic, a contradiction with Lemma~\ref{cyclic} by the definition of $y$.
		\par

		Suppose now $n\geq 2$ and that $\lambda_{1}\neq 0$. If $\lambda_{i}\neq 0$ for any $i>1$ then replace $a_i$ by $\lambda_{1}a_i-\lambda_{i}a_1$. So, we may suppose that $\lambda_{i}=0$ for $i>1$. Put $y_1=a_1+x^2$, $y_2=b+x$ where $b=a_i$ for some $i>1$. Then $\langle y_1\rangle=Fy_1$, $\langle y_2\rangle=\sum_{i=1}^kF(b+x^i)$. Now we must have
		\[ [y_1,y_2]=a_1+\lambda_{1}x^k+x^3=a_1+x^2+\sum_{i=1}^k\mu_i(b+x^i).
		\] Equating coefficients we get:
		\par
		
		\noindent { if $k>3$,} then $\mu_1=0$, $\mu_2=-1$, $\mu_3=1$, $\mu_4=\cdots=\mu_{k-1}=0$, $\mu_k=\lambda_{1}$ and $\sum_{i=1}^k \mu_i=0$, whence $\lambda_{1}=0$, a contradiction;
		\par
		
		\noindent { if $k=3$,} then $\mu_1=0$, $\mu_2=-1$, $\mu_3=1+\lambda_{1}$  and $\sum_{i=1}^3 \mu_i=0$, whence $\lambda_{1}=0$, a contradiction;
		\par
		
		\noindent { if $k=2$,} then $1+\lambda_{1}=1$, whence $\lambda_{1}=0$, a contradiction.
		\par
		
		Then, the claim is proved and $L$ is as in (ii).

                    Next we consider the sufficiency. Suppose first that (i) holds, so $L=A+Fx$ is a non-Lie almost abelian Leibniz algebra. Let $y_1=a_1+\lambda x, y_2=a_2+\mu x\in L$. Then
\[  [y_1,y_2]=\mu a_1=\left\{ \begin{matrix}  \mu(a_1+\lambda x) & \hbox{if} & \lambda=0 \\ \mu\lambda^{-1}(a_1+\lambda x)^2 & \hbox{if}&\lambda\neq 0 \end{matrix}\right\}\in \langle y_1\rangle. 
\] It follows from Lemma \ref{qi} that every subalgebra is a weak quasi-ideal.
\par

Finally, suppose that (ii) holds. We have that $\phi(L)=\sum_{i=2}^kFx^i$ if $A\neq 0$ and $\phi(L)=\sum_{i=2}^{k-1} (x^i-x^{i+1})$ if $A=0$. In either case, $L/\phi(L)$ is a non-Lie almost abelian Leibniz algebra. Let us prove now that every subalgebra of $L$ is a weak quasi-ideal. Suppose $A\neq 0$; otherwise, Lemma~\ref{cyclic} yields the result. Take $y=a+\sum_{i=1}^k\lambda_ix^i$ and $z=b+\sum_{i=1}^k\mu_ix^i$ for some $a,b\in A$; then \[[y,z]=\mu_1 a + \mu_1\sum_{i=1}^{k}\lambda_i x^{i+1}.\] Consider first the case $\lambda_1\neq 0$; then, $[x,y]=\mu_1 \lambda_1^{-1} y^2\in \langle y\rangle$. On the contrary, if $\lambda_1=0$, we have $\langle y\rangle=F y$. Take $\mu_1=1$ without loss of generality. 
Then $z^r=b+\sum_{i=1}^{k}\mu_i x^{i+r-1}$ for all $r\geq 1$. Set $\mu_0=0$ and define $w_r=z^r-z^{r+1}=x^r+\sum_{i=1}^{k-1-r}(\mu_{i+1}-\mu_i)x^{r+1}-\mu_{k-r}x^k$ for $1\leq r \leq k$. Also, we define recursively
\begin{align*}
	\gamma_2&=-\lambda_2; \\
	\gamma_r&=\lambda_{r-1}-\lambda_r-\sum_{i=2}^{r-1}\gamma_i(\mu_{r+1-i}-\mu_{r-i})\qquad \text{for } 3\leq r\leq k.
\end{align*}
We claim that \[[y,z]=y + \sum_{i=2}^{k}\gamma_i w_i\in \langle y\rangle + \langle z\rangle;\]
indeed, one can see that, by the choice of the $\gamma_i$, it suffices to check that $\sum_{i=2}^{k-1}\gamma_i\mu_{k-i}=-\lambda_{k-1}$. This can be done by an easy induction in $k$. By Lemma~\ref{qi}, every subalgebra is a weak quasi-ideal.


	\end{proof}
	
	\begin{theor}\label{sqrt} Let $L$ be a Leibniz algebra over a field $F$ of characteristic different from $2$ and with $\sqrt{F}\subseteq F$. Then every subalgebra is a weak quasi-ideal and $L/\phi(L)$ is an almost abelian Lie algebra if and only if $L= A\dot{+} \langle x\rangle$ where $A\neq 0$ is an abelian ideal of $L$, $\langle x\rangle$ is a nilpotent cyclic Leibnix algebra, $[a,x]=a$ and $[x,a]=-a$, $[x^k,a]=[a,x^k]=0$ for all $a\in A, k>1$. 
	\end{theor}
	\begin{proof} Let $L/\phi(L)$ be an almost abelian Lie algebra. Then $I\subseteq \phi(L)$. But $L/I$  is an abelian or almost abelian Lie algebra, by \cite[Theorem 3.8]{amayo}, so $\phi(L)\subseteq I$, whence $\phi(L)=I$. We proceed, as in Theorem \ref{non-Lie}, by induction on $\dim \phi(L)$. The result clearly holds if $\dim \phi(L)=0$, with $\langle x \rangle = Fx$, so suppose it holds whenever $\dim \phi(L)\leq  m$ ($m\geq 0$) and let $L$ be such that $\dim \phi(L)=m+1$. Since $L$ is supersolvable, there is a one-dimensional ideal $Fa_n$ inside $I$, and $L/Fa_n$ is of the given form.
		\par
		
	 Then $L$ will have a basis $a_1, \ldots, a_{n},$ $x,\ldots,x^k$ where $[a_i,x]=a_i+\alpha_i a_n$, $[x,a_i]=-a_i+\beta_i a_n$  for all $1\leq i\leq n-1$, $[a_i,a_j]=\gamma_{ij}a_n$ for $1\leq i,j \leq n$, $x^{k+1}=\delta a_n$ and $[a_n,x]=\epsilon a_n$.  Now, $[x,[x,a_i]] =-[x,a_i]$ for $1\leq i\leq n-1$ since $a_n\in I$, so, by replacing $a_i$ by $[x,a_i]$ we may assume that $\beta_i=0$. Also, $$-a_i=[x,a_i]=[x,[a_i,x]]=[[x,a_i],x]-[x^2,a_i]=-a_i-\alpha_i a_n-[x^2,a_i],$$ so
		\begin{align*}
		[x^2,a_i]=-\alpha_i a_n \hbox{ for } 1\leq i\leq n-1. 	
		\end{align*}
 A straightforward induction argument shows that
 \begin{align}\label{alpha}
[x^{ r},a_i]=(-1)^{r+1}\alpha_i(1-\epsilon)^{r-2} a_n \hbox{ for }   r>1.
 \end{align}
\par

Also consider
\begin{align*}\gamma_{ni}a_n=&[a_n,a_i]=[a_n,[a_i,x]]=[[a_n,a_i],x]-[[a_n,x],a_i]\\=&\gamma_{ni}[a_n,x]-\epsilon[a_n,a_i]=0,\end{align*}
so $\gamma_{ni}=0$ for all $1\leq i\leq n-1$.
Now
		
		\[\epsilon\gamma_{ij}a_n=[[a_i,a_j],x]=[a_i,[a_j,x]]+[[a_i,x],a_j]=2\gamma_{ij} a_n,\] so $\epsilon=2$ or $\gamma_{ij}=0$ for all $1\leq i,j\leq n-1$.
		
		Suppose that $\delta=0$, so $a_n \notin \langle x\rangle$ and Equation~\eqref{alpha} implies that  $\epsilon=1$ or $\alpha_i=0$ for $1\leq i\leq n-1$. If $\epsilon=1$, then $\gamma_{ij}=0$ for $1\leq i\leq n$. But now $\sum_{i=1}^{n-1}Fa_i+\langle x\rangle$ is a maximal subalgebra, so $a_n \notin \phi(L)$, a contradiction. So suppose $\alpha_i=0$ for $1\leq i\leq n-1$. 
		
		We claim that, if $n\geq 3$, then there exists an $a\in A$ such that $a^2=0$. Indeed, if this was not the case, consider \[(\lambda a_1+\mu a_2)^2=(\lambda^2\gamma_{11}+\mu^2\gamma_{22}+\lambda\mu(\gamma_{12}+\gamma_{21}))a_n,\] with $\gamma_{11},\gamma_{22}\neq 0$. Since $\sqrt{F}\subseteq F$ we can find $a=\lambda a_1+\mu a_2$ such that $a^2=0$. 
Then  $[a+a_n,x]=a+\epsilon a_n$ implies that $\epsilon =1$, and we get a contradiction as before.

		If $n= 2$ and  $\gamma_{11}= 0$, we argue as in the previous paragraph. Otherwise, $\epsilon=2$ and put $a_2=a_1^2$. Take $y_1=a_1-\frac{1}{2} a_2 + x$ and $y_2=a_1-\frac{1}{4} a_2 + 2x$. Elementary calculations show that $[y_1,y_2]=a_1-a_2+2x^2\notin \langle y_1\rangle + \langle y_2\rangle=F y_1 + F y_2 +\sum_{i=2}^k F x^i$, so by Lemma~\ref{qi} $\langle y_1,y_2\rangle $ is not a weak quasi-ideal, a contradiction.
		\par

		So suppose now that $\delta\neq 0$. Then we can put $a_n=x^{k+1}$ and $x^{k+2}=0$ or $x^{k+2}=x^{k+1}$. In either case, $\epsilon\neq 2$, so $\gamma_{ij}=0$ for $1\leq i,j\leq n$  as before.
		\par

		Suppose that $\epsilon\neq 0$. Then, $[a_i,x-x^2]=a_i+\alpha_ix^{k+1}$, so $\alpha_ix^{k+1}\in \langle x - x^2\rangle=F(x-x^2)+\dots+F(x^k-x^{k+1})$. It follows that $\alpha_i=0$ for $1\leq i \leq n-1$ and  $\sum_{i=1}^{n-1}Fa_i+\langle x-x^2\rangle$ is a maximal subalgebra.  But, again, $a_n \notin \phi(L)$, so this case cannot arise.
		\par

		So, we must have that $a_n=x^{k+1}$ and $x^{k+2}=0$. Now (1) shows that $\alpha_i=0$ for $1\leq i\leq n-1$, and, putting $A=Fa_1+\cdots+Fa_{n-1}$, we see that $L$ is of the given form.
\par

It remains to show sufficiency.  Consider $y_1=a_1+\sum_{i=1}^{k}\lambda_ix^i$, $y_2=a_2+\sum_{i=1}^{k}\mu_ix^i$. If $\lambda_1\neq 0$ or $\mu_1\neq 0$, a straightforward calculation shows that $\langle y_1\rangle + \langle y_2 \rangle =F(a_1+\lambda_1 x) + F (a_2 + \mu_1 x) + \sum_{i=2}^k Fx^i$. Then, \begin{align*}[y_1,y_2]&=\mu_1a_1-\lambda_1a_2+\sum_{i=1}^{k-1}\lambda_i\mu_1x^{i+1}\\ &=\mu_1(a_1+\lambda_1x)-\lambda_1(a_2+\mu_1x)+\sum_{i=1}^{k-1}\lambda_i\mu_1x^{i+1}\in \langle y_1\rangle +\langle y_2\rangle,
\end{align*}
and every subalgebra of $L$ is a weak quasi-ideal by Lemma~\ref{qi}. Also $\phi(L)= \sum_{i=2}^k Fx^i$, so $L/\phi(L)$ is a Lie almost abelian algebra.
	\end{proof}

\begin{remark} \emph{
Note that in~\cite[Proposition 6.3]{st} it was proven that, for a Leibniz algebra $L$ over an algebraically closed field, the following are equivalent:
\begin{enumerate}
\item[(i)] $L$ is modular;
\item[(ii)] every subalgebra of $L$ is a weak quasi-ideal of $L$.
\end{enumerate}
It is possible to extend this equivalency to:
\begin{enumerate}
	\item[(iii)] $L$ is upper semi-modular.
\end{enumerate}
Indeed, suppose that $L$ is upper semi-modular. By~\cite[Corollary 4.5 and Proposition 5.1]{st}, $L$ is supersolvable and hence lower semi-modular. The general theory of lattices of finite length says that, being upper and lower semi-modular, then $L$ is also modular, see~\cite{schmidt}.}
\end{remark}

If $L$ is solvable then the restriction on $F$ is unnecessary in the above remarks.

\section{Modular Symmetric Leibniz algebras}
Finally, we have the following classification result for modular symmetric Leibniz algebras.

\begin{theor} Let $L$ be a symmetric Leibniz algebra over an algebraically closed field of characteristic not two. Then $L$ is modular if and only if one of the following condition holds:
	\item[(i)] $L$ is an abelian Lie algebra;
	\item[(ii)] $L$ is an almost abelian Lie algebra;
	\item[(iii)]  $L=E\oplus Z$, where $Z$ is a central ideal of $L$ and  $E$ is an extraspecial Leibniz algebra such that $J=\{x\in L \mid x^2=0\}$ is an abelian ideal of $L$;
	\item[(iv)] $L=B+F y + F y^2$, where $B$ is an abelian ideal, $[b,y]=b=-[y,b]$ for all $b\in B$, and $Z(L)=Fy^2$.
\end{theor}

\begin{proof} We first prove the necessity part. Suppose that $L$ is modular. Since $L$ is symmetric, note that $I$ is a central ideal of $L$. Moreover, $L/I$ is a modular Lie algebra and then, by  \cite[Corollary 2.2]{kol}, $L/I$ is abelian or almost abelian. In the former case, $L$ is nilpotent of class at most 2, and it follows from \cite{st} that $L$ is either abelian or of the form given in condition (iii).
	
	Suppose next that $L/I$ is almost abelian. Then $L=F y + A$, where $y\in L$, $A$ is an ideal of $L$ with $I\subseteq A$, $A/I$ is abelian and $[a,y]=a + i_a$ for all $a\in A$, with $i_a\in I$.  As $A$ is nilpotent of class at most 2, by \cite{st} we easily deduce that $I= F y^2 + F a^2$ and $[A,A]= F a^2$ for some $a\in A$. For every $x\in A$ we have 
	$$[x,y]+[y,x]=[[x,y], y ]-[ y,[y,x]]= -[y,[x,y]]+[[y,x],y]=0,$$ 
	hence $[x,y]=-[y,x]$. 
	
	We show that $A$ is actually abelian. 
	Suppose the contrary, so that $a^2\neq 0$. 
	We have $[a,y]=a+\lambda a^2+\mu y^2$ for some $\lambda, \mu \in F$, so we can replace $a$ by $a+\lambda a^2+\mu y^2$ to assume that $[a,y]=a$. But then, as Leibniz elements in symmetric Leibniz algebras are central, we have $0=[a^2,y]=2a^2$. Since $F$ has characteristic different from 2, we conclude that $a^2=0$, a contradiction. 
	
	Now, if $y^2=0$ then $L$ is an almost abelian Lie algebra.   
	Suppose then $y^2\neq 0$. For every $a\in A$ we have $[a,y]=a+\mu y^2$ for some $\mu \in F$. Therefore, we can find $b_1,b_2,\ldots,b_n \in A$ such that $[b_i,y]=b_i$ and $[b_i,b_j]=0$ for all $i,j$ and $  y^2,b_1,\ldots,b_n$ is a basis of $A$. Then $B=F b_1 + \cdots + F b_n$ is an abelian ideal of $L$ and condition (iv) of the statement holds.
	
	The proof of the sufficiency part follows from  \cite[Proposition 1.2]{kol}, \cite[Proposition 6.7]{st}, Lemma \ref{qi}, and Theorem \ref{sqrt}. 
\end{proof}

\section*{Acknowledgments}

The first author was supported by the Austrian Science Foundation FWF, grant P 33811, while working on this project.


\begin{thebibliography}{1}
\bibitem{amayo} {\sc R.K.Amayo}, `Quasi-ideals of Lie algebras I', {\em Proc. London Math. Soc. (3)} {\bf 33} (1976), 28--36.

\bibitem{1} {\sc R.K. Amayo and J. Schwarz}, `Modularity in Lie Algebras', {\em Hiroshima Math. J.} {\bf 10} (1980), 311-322.
		
\bibitem{barnes} {\sc D.W. Barnes}, `Some theorems on Leibniz algebras', {\em Comm. Algebra} {\bf 39(7)} (2011), 2463--2472.

\bibitem{bloh} {\sc A. Bloh}. `On a generalization of the concept of Lie algebra'. {\em Dokl. Akad. Nauk SSSR.} {\bf 165} (1965), 471--473.

\bibitem{7} {\sc K. Bowman and D.A.Towers}, `Modularity conditions in Lie algebras', {\em Hiroshima Math. J.} {\bf 19} (1989), 333-346.

\bibitem{14} {\sc K. Bowman and V.R. Varea}, `Modularity$^\ast$ in Lie algebras', {\em Proc. Edin. Math. Soc.} {\bf 40(2)} (1997), 99-110.

\bibitem{2} {\sc A.G. Gein}, `Modular rule and relative complements in the lattice of subalgebras of a Lie algebra', {\em Sov. Math.} {\bf 31(3)} (1987), 22-32; translated from {\em Izv. Vyssh. Uchebn. Zaved. Mat.} {\bf 83} (1987), 18-25.

\bibitem{8} {\sc A.G. Gein}, `Semimodular Lie algebras', {\em Siberian Math. J.} {\bf 17} (1976), 243-248; translated from {\em Sibirsk Mat. Z.} {\bf 17} (1976), 243-248.

\bibitem{3} {\sc A.G. Gein}, `On modular subalgebras of Lie algebras', {\em Ural Gos. Univ. Mat. Zap.} {\bf 14} (1987), 27-33.

\bibitem{10} {\sc A.G. Gein and V.R. Varea}, `Solvable Lie algebras and their subalgebra lattices', {\em Comm.  Alg.} {\bf 20(8)} (1992), 2203-2217. Corrigenda: {\em Comm. Alg.} {\bf 23(1)} (1995), 399-403.

\bibitem{kol} {\sc B. Kolman}, 'Semi-modular Lie algebras', {\em J. Sci. Hiroshima Univ. Ser. A-I} {\bf 29} (1965), 149-163.

\bibitem{11} {\sc A.A. Lashi}, `On Lie algebras with modular lattices of subalgebras', {\em J. Algebra} {\bf 99} (1986), 80-88.

\bibitem{loday1} {\sc J.-L. Loday}, `Une version non commutative des alg\`{e}bres de Lie: les algèbres de Leibniz'. {\em Enseign. Math. (2)} {\bf 39 (3–4)} (1993), 269--293.

\bibitem{loday2} {\sc J.-L. Loday and  T. Pirashvili}, `Universal enveloping algebras of Leibniz algebras and (co)homology', {\em Math. Annalen} {\bf 296 (1)} (1993) 139--158.

\bibitem{sch} {\sc C. Scheiderer}, `Intersections of maximal subalgebras in Lie algebras', {\em J. Algebra} {\bf 105} (1987), 268-270.
		
\bibitem{schmidt}  {\sc R. Schmidt}, `Subgroup Lattices of Groups', de Gruyter, Berlin, New York (1994).
		
\bibitem{st} {\sc S. Siciliano and D.A. Towers}, `On the subalgebra lattice of a Leibniz algebra', {\em Comm. Alg.} {\bf 50 (1)} (2021), 255--267.
		
\bibitem{frat} {\sc D.A. Towers}, `A Frattini theory for algebras', {\em
			Proc. London Math. Soc.} (3) {\bf 27} (1973), 440--462.

\bibitem{13} {\sc D.A. Towers}, `Semimodular subalgebras of a Lie algebra', {\em J. Algebra} {\bf 103} (1986), 202-207.

\bibitem{15} {\sc D.A. Towers}, `On modular$^\ast$ subalgebras of a Lie algebra', {\em J. Algebra} {\bf 190} (1997), 461-473.

\bibitem{4} {\sc V.R. Varea}, `Modular subalgebras, quasi-ideals and inner ideals in Lie Algebras of prime characteristic', {\em Comm. Alg.} {\bf 21(11)} (1993), 4195--4218.

\bibitem{5} {\sc V.R. Varea}, `The subalgebra lattice of a supersolvable Lie algebra', {\em In Lecture Notes in Mathematics}. Springer-Verlag:
New York, 1989; Vol. 1373, 81-92.

\bibitem{6} {\sc V.R. Varea}, `On modular subalgebras in Lie algebras of prime characteristic, {\em Contemporary Math.} {\bf 110} (1990), 289-307.

\bibitem{9} {\sc V.R. Varea}, `Lower Semimodular Lie algebras', {\em Proc. Edin. Math. Soc.} {\bf 42}(1999), 521-540.

\bibitem{12} {\sc V.R. Varea}, `Lie algebras whose maximal subalgebras are modular', {\em Proc. Roy. Soc. Edinburgh Sect. A} {\bf 94} (1983),
9-13.
		
		
	\end{thebibliography}
\end{document}